\documentclass[leqno,11pt]{article}
\usepackage{amsmath, amsthm}
\usepackage{amsfonts}
\usepackage{amssymb}
\usepackage{graphicx}
\usepackage{color}
\usepackage{hyperref}

\newtheorem{theorem}{Theorem}[section]
\newtheorem{corollary}{Corollary}
\newtheorem{lemma}[theorem]{Lemma}

\newcommand{\E}{{\mathcal E}}
\newcommand{\re}{\mathbb{R}}
\newcommand{\ren}{\mathbb{R}^N}

\newcommand{\dint}{\displaystyle\int}
\newcommand{\dx}{\,{\rm d}x}

\newcommand{\dt}{\,{\rm d}t}
\newcommand{\rd}{{\rm d}}
\newcommand{\A}{\mathcal{L}}

\def\qed{\,\unskip\kern 6pt \penalty 500
\raise -2pt\hbox{\vrule \vbox to8pt{\hrule width 6pt
\vfill\hrule}\vrule}\par}

\title{\bf Recent progress in the theory of Nonlinear \\Diffusion with Fractional Laplacian Operators}
\author{ {\large Juan Luis V\'azquez} \\[10pt]
Universidad Aut\'onoma de Madrid, Spain}

\date{November 2013}
\begin{document}
\maketitle

\bigskip

\begin{abstract}
We report on recent progress in the study of nonlinear diffusion equations involving nonlocal, long-range diffusion effects. Our main concern is the so-called fractional porous medium equation, $\partial_t u +(-\Delta)^{s}(u^m)=0$, and some of its generalizations. Contrary to usual porous medium flows, the fractional version has infinite speed of propagation for all  exponents $0<s<1$ and $m>0$; on the other hand, it also generates an $L^1$-contraction semigroup which depends continuously on the exponent of fractional differentiation and the exponent of the nonlinearity.

 After establishing the general existence and uniqueness theory, the main properties are described: positivity, regularity, continuous dependence, a priori estimates, Schwarz symmetrization, among others. Self-similar solutions are constructed (fractional Barenblatt solutions) and they explain the asymptotic behaviour of a large class of solutions. In the fast diffusion range we study extinction in finite time and we find suitable special solutions. We discuss KPP type propagation.  We also examine some related equations that extend the model and briefly comment on current work.
 \end{abstract}

\section{Introduction}\label{sec.intro}

This is a follow-up to a previous survey by the author \cite{VazAbel}, which reported on recent research on two models of nonlinear diffusion processes involving long-range diffusion, written on the occasion of the Abel Symposium held in Oslo in 2010. Such evolution processes are  represented by nonlinear parabolic equations involving nonlocal operators of the fractional Laplacian type.

\medskip

Recapitulating what was said there, the classical theory of diffusion is expressed mathematically by means of the heat equation, and more generally by parabolic  equations of linear  type; it  has  had an enormous success and is now a foundation stone in science and technology. The last half of the past century has witnessed intense activity and progress in the theories of nonlinear diffusion, examples being the Stefan Problem, the Porous Medium Equation, the $p$-Laplacian equation, the  Total Variation Flow, evolution problems of Hele-Shaw type, the Keller-Segel chemotaxis system, and many others. Reaction diffusion has also attracted considerable attention. In the last decade there has been a surge of activity focused on the use of so-called fractional diffusion operators to replace the standard Laplace operator (and other kinds of elliptic operators with variable coefficients), with the aim of further extending the theory by taking into account the presence of the long range interactions that occur in a number applications. The new operators do not act by pointwise differentiation but by a global integration with respect to a singular kernel;  in that way the nonlocal character of the process is expressed.

\subsection{Fractional operators}\label{subs.frop}

Though there is a wide class of interesting nonlocal operators under scrutiny, a substantial part of the current work deals diffusion modeled by the so-called fractional Laplacians. We recall that the fractional Laplacian operator is a kind of isotropic differentiation operator of order $2s$, for some $s\in (0,1)$, that can be conveniently defined  through its Fourier Transform symbol, which is $|\xi|^{2s}$.  Thus, if $g$ is a function in the Schwartz class in $\ren$, $N\ge 1$, we write  $(-\Delta)^{s}  g=h$ if
\begin{equation}
\widehat{h}\,(\xi)=|\xi|^{2s} \,\widehat{g}(\xi)\,,
\label{def-fourier}
\end{equation}
so that for $s=1$ we recover the standard Laplacian. This definition allows for a wider range of parameters $s$. The interval of interest for fractional diffusion is $0<s\le 1$, and for $s< 1$  we can also use the integral representation
\begin{equation}\label{def-riesz}
(-\Delta)^{s}  g(x)= C_{N,2s }\mbox{
P.V.}\int_{\mathbb{R}^N} \frac{g(x)-g(z)}{|x-z|^{N+{2s}
}}\,dz,
\end{equation}
where P.V. stands for principal value and $C_{N,\sigma}$ is a normalization constant, with precise value $C_{N,2s}=2^{2s}s\Gamma((N+2s)/2)/( \pi^{N/2}\Gamma(1-s))$. In the limits $s\to 0$ and $s\to 1$ it is possible to recover respectively the identity or the standard minus Laplacian, $-\Delta$, cf. \cite{Brezis, Mazya}. Remarkably, the latter one cannot be represented by a nonlocal formula of the type \eqref{def-riesz}.
It is also useful to recall that the operators $(-\Delta)^{-s}$, $0< s<1$,  inverse of the former ones,  are  given by standard convolution expressions:
\begin{equation}
(-\Delta)^{-s}  g(x)= C_{N,-2s }\int_{\mathbb{R}^N}
\frac{g(z)}{|x-z|^{N-{2s} }}\,dz, \label{def-riesz2}
\end{equation}
in terms of the usual Riesz potentials. Basic references for these operators are the books by
Landkof \cite{Landkof} and Stein \cite{Stein}. A word of caution: in the literature we often find the notation $\sigma=2s$, and then the desired interval is $0<\sigma<2$. According to that practice, we will sometimes use $\sigma$ instead of $s$.

\smallskip

 Another option is to use the classical way of defining the fractional powers of a linear self-adjoint nonnegative operator, and it is expressed in terms of the  semigroup associated to such an operator. In the case of the standard Laplacian operator, it reads
\begin{equation}\label{sLapl.Rd.Semigroup}
\displaystyle(-\Delta)^{s}
g(x)=\frac1{\Gamma(-s)}\int_0^\infty
\left(e^{t\Delta}g(x)-g(x)\right)\frac{dt}{t^{1+s}}.
\end{equation}
All the above definitions are equivalent when dealing with the Laplacian on the whole space $\re^N$.

\smallskip

The interest in these fractional operators has a long history in Probability for reasons we explain in the next paragraph. Motivation from Mechanics appears in the famous Signorini problem (with $\alpha=1/2$), cf. \cite{Signor, CaffSign}. And there are applications in Fluid Mechanics, cf. \cite{CaffVass, KisNazVol}. There is a wide literature on the subject, both for its relevance to Analysis, PDEs, Potential Theory, Stochastic Processes, and Finance, and for  the growing number of practical applications.  See e.g. \cite{Caffarelli-Silvestre, DPV11, Valdinoc, VazAbel} where further references can be found.

\smallskip

The systematic study of the corresponding PDE models with fractional operators is relatively recent, and many of the results have been established in the last decade. A part of the current research concerns linear or quasilinear equations of elliptic type. This is a huge subject with well-known classical references that will not be discussed here.


\subsection{\bf Linear evolution processes}

A significant  part of the motivation and also a large part of the recent literature is related to evolution problems, that we will discuss in the sequel. Thus, the difference between the standard and the fractional Laplacian is best seen in the stochastic point of view (the stochastic theory of diffusion), and it consists in re-examining the conditions under which the Brownian motion is derived, and taking into account long-range interactions instead of the usual interaction driven by close neighbors. This change of  model explains characteristic new features of great importance, like enhanced propagation with the appearance of fat tails at long distances (such tails are to be compared with the typical exponentially small tails of the standard diffusion, or the compactly supported solutions of porous medium flows). Moreover, the space scale of the propagation of the distribution is not proportional to $t^{1/2}$ as in the Brownian motion, but to another power of time, that can be adjusted in the model; this is known as {\sl anomalous diffusion.} The fractional Laplacian operators of the form $(-\Delta)^{\sigma/2}$, $\sigma\in(0,2)$, are actually the infinitesimal generators of stable L\'{e}vy processes \cite{Applebaum, Bertoin, CKS2010}.

\smallskip

We will be concerned with evolution partial differential equations that combine diffusion and fractional operators. As already mentioned, a  great variety of diffusive problems in nature, namely those referred to as normal diffusion, are satisfactorily described by the classical Heat Equation or Fokker-Planck linear equation. However, anomalous diffusion is nowadays intensively studied, both theoretically and experimentally since it conveniently explains a number of phenomena in several areas of physics, finance, biology, ecology, geophysics, and many others, which can be  briefly summarized as having non-Brownian scaling. This leads to linear anomalous diffusion equations, see e.g. \cite{Applebaum, ContTankov2004, Woy2001}. Fractional kinetic equations of the diffusion, diffusion-advection, and Fokker-Planck type represent a useful approach for the description of transport dynamics in complex systems which are governed by anomalous diffusion. These fractional equations are usually derived asymptotically from basic random walk models, cf.  \cite{Jara0, JKOlla, MK2000, MMM,Valdinoc, VIKH, WZ} and their references.

\smallskip

 The standard linear evolution equation involving fractional diffusion is
\begin{equation}\label{linfraceq}
\frac{\partial u}{\partial t}+(-\Delta)^{s}(u)=0\,,
\end{equation}
a usual model for anomalous diffusion.  The equation is solved with the aid of well-known functional analysis tools, like Fourier transform; for instance, it is proved that it generates a semigroup of ordered contractions in $L^1(\mathbb{R}^N)$.  Moreover, in this setting it has the integral representation
\begin{equation}
u(x,t)=\int_{\mathbb{R}^N}K_s(x-z,t)f(z)\,dz\,,
\label{lineal}
\end{equation}
where $K_s$ has Fourier transform  $\widehat K_s(\xi,t)=e^{-|\xi|^{2s} t}$. This means that, for $0<s<1$, the kernel $K_s$ has the form
\begin{equation}
K_s(x,t)=t^{-N/2s}F(|x|/t^{1/2s})
\end{equation}
 for some profile function $F=F_s$ that is positive and decreasing, and it behaves at infinity
like $F(r)\sim r^{-(N+2s)}$, \cite{Blumenthal-Getoor}. When $s=1/2$, $F$ is explicit:
\begin{equation}\label{kernel.lin}
F_{1/2}(r)={C}\,(a^2+r^2)^{-(N+1)/2}\,.
\end{equation}
If $s=1$ the function $K_{s=1}$ is the Gaussian heat kernel, which has a negative square exponential tail, i.e., a completely different asymptotic behavior.

\smallskip

An integral representation of the evolution of the form \eqref{lineal} is not available in the nonlinear models coming from the applications. This implies the need for new methods, thus motivating our work to be described below.

\section{Nonlinear diffusion models}\label{sec.nd}

A main feature of current research in the area of PDEs  is the interest in nonlinear equations and systems. There are a number of models of  evolution equations that can be considered nonlinear counterparts of the linear fractional heat equation, and combine Laplace operators and nonlinearities in different ways. Let us mention some of the most popular in the recent PDE literature.

\smallskip

\noindent $\bullet$ {\bf Type I.} A natural option is to consider the equation
\begin{equation}\label{fpme}
\partial_t u +(-\Delta)^{s}(u^m)=0
\end{equation}
with $0<s<1$ and $m>0$. This is mathematically the fractional version of the standard Porous Medium Equation (PME)
\begin{equation}
\partial_t u =\Delta(u^m)\,,
\end{equation}
that is recovered as the limit $s\to 1$ and has been extensively studied, cf. \cite{ArBk86, Vapme}. We will call equation \eqref{fpme} the {\sl Fractional Porous Medium Equation, FPME,} as proposed in our first works on the subject \cite{pqrv1, pqrv2}, in collaboration with de A. de Pablo,  F. Quir\'os and A. Rodr\'{\i}guez.

\smallskip

Interest in studying the nonlinear model we propose is two-fold: on the one hand, experts in the mathematics of diffusion want to understand the combination of fractional operators with porous medium type propagation (which is described as degenerate parabolic). Models of this kind arise in statistical mechanics when modeling for instance heat conduction with anomalous properties. On the other hand, it is mentioned in heat control problems by \cite{ACld}. The rigorous study of such nonlinear models has been delayed until this time by the mathematical difficulties in treating at the same time the nonlinearity and fractional diffusion.

\smallskip

We will devote most of this paper to explain the main results obtained so far on the mathematical theory of this equation.  Thus, Section \ref{sec.fpme} contains the theory of existence, uniqueness and continuous dependence of the Cauchy problem posed in the whole space $\ren$, developed in \cite{pqrv1, pqrv2}.  A main result is the property of infinite speed of propagation. Let us explain the result at this moment and point out the contrast to the PME. While one of the best known features of the PME for $m>1$ is the fact that compactly supported nonnegative initial data give rise to solutions that are also compactly supported as functions of $x$ for every fixed $t>0$, this property does not hold in equation \eqref{fpme} if $s<1$. Indeed, the unique nonnegative solutions we have constructed for nontrivial data $u(x,0)\in L^1_+(\ren)$ are positive everywhere for all $t>0$. This happens for all $m>0$ and all $0<s<1$ (for $m$ close to 0 attention must be paid to extinction of the solution in finite time, but that is another issue, see below). Some traits are common with the PME equation: an $L^1$-contraction semigroup is constructed and it depends continuously on the exponent of fractional derivation and the exponent of the nonlinearity.

\smallskip

The further investigation of the FPME has taken a number of directions. Thus, in Section \ref{sec.bar} we discuss the existence of self-similar solutions with conserved finite mass (we call them {\sl fractional Barenblatt solutions}), following the study made  in \cite{VazBar2012}. The main  properties are derived, allowing to construct the self-similar profiles. Generally, such profiles are not explicit. These solutions are then used to explain the asymptotic behavior of the whole class of nonnegative solutions with finite mass (a kind of fractional central limit theorem).

\smallskip

The study of a priori estimates (some of them universal), quantitative bounds on positivity and Harnack estimates is done in Section \ref{sec.mb} following the results obtained in \cite{BV2012} in collaboration with M. Bonforte. This is done in the favorite setup, the Cauchy problem. But the fractional Laplace operators offer many novelties when posed in a bounded domain of the Euclidean space, say with Dirichlet boundary conditions. We report in Section \ref{sec.dir} on recent work done with M. Bonforte on the issue and comment on related works.

\smallskip

As a natural extension of the equation, we have recently considered equations of the form
\begin{equation}\label{gfpme}
\partial_t u +(-\Delta)^{s}\Phi(u)=0
\end{equation}
where $\Phi$ is a monotone increasing function. Collaboration with Pablo, Quir\'os and Rodr\'{\i}guez, \cite{pqrv3} and  \cite{vpqr}; regularity is a main issue, we review the results in Section \ref{sec.gnl}. Section \ref{sec.symm} reports on the technique of Schwarz symmetrization of the solutions of the these problems and derives a priori estimates. This is the result of collaboration with B. Volzone \cite{VazVol, VazVol2}. Finally, numerical methods for these equations are being investigated by our team and a number of authors, and we will make a comment on that issue.

\smallskip

\noindent $\bullet$  The combination of fractional diffusion with other effects leads to interesting behaviour. One of the most common combinations is reaction-diffusion. In that sense there has  been interesting work on the possible extension of the well-known KPP behaviour to cover fractional diffusion, both in the linear and nonlinear case. Results by X. Cabr\'e and J. M. Roquejoffre \cite{CabreRoquejoffre2} (for linear diffusion) and Diana Stan and the author \cite{StanVazquezKPP} (in the nonlinear case) are covered in Section \ref{sec.kpp}.

\smallskip

On the other hand, the combination of fractional diffusion convection appears in a number of models of the so-called geostrophic equations, like
\begin{equation}
\partial_t u + {\bf v}\cdot \nabla u =(-\Delta)^s u,
\end{equation}
where $0<s<1$ and $\bf v(x,t)$ is a divergence free vector field related to $u$ in different ways. Important works in this direction related to our research are due to Caffarelli and Vasseur \cite{CaffVass} and Kiselev et al. \cite{KisNazVol}, but again, we will not deal with this topic here.  This last model uses linear fractional diffusion. Nonlinear diffusion combined with convection has been studied by some authors, like Cifani and Jakobsen \cite{CJ11}, where references to other models are given.

\smallskip

\noindent $\bullet$  {\bf Type II.} We recall here the other model of diffusion with fractional operators where the  author has been involved, mainly in collaboration with Luis Caffarelli. Since it was extensively reported in the survey paper \cite{VazAbel} we will only give a short reminder. This alternative model is derived in a more classical way from the Porous Medium Equation since it is based on the usual Darcy  law, with the novelty that the pressure is related to the density by an inverse fractional Laplacian operator. This nonlinear fractional diffusion equation of porous medium type takes the form
\begin{equation}\label{PMEfp}
u_t=\nabla (u\, \nabla {\mathcal K} u),
\end{equation}
where $\mathcal K$ is the Riesz operator that typically expresses the inverse to the fractional Laplacian, ${\mathcal K} u=(-\Delta)^{-s}u$. This has been studied by Caffarelli and the author in \cite{CV1, CV2} and Biler, Karch and Monneau \cite{BKM, BIK}, where the equation is derived in the framework of the theory of dislocations. We recall that this model has some strikingly different properties, like lack of strict positivity and occurrence of free boundaries. Self-similar solutions also exist but their existence and properties are quite different from those of equation \eqref{fpme}. Even the asymptotic behaviour is quite different. In order to distinguish both models the name {\sl porous medium equation with fractional pressure} has been proposed for equation \eqref{PMEfp}.

\smallskip

As for recent work, the boundedness and H\"older regularity of finite mass solutions is studied in joint work  with  L. Caffarelli and F. Soria, \cite{CSV}.  As a limit case of this second model $s\to1$,  one  obtains a variant of the equation for the evolution of vortices in superconductivity  derived heuristically  by Chapman-Rubinstein-Schatzman \cite{CRS} and W. E \cite{WE} as the hydrodynamic limit of Ginzburg Landau,  and studied by  Lin  and Zhang \cite{Liz}, and Ambrosio and Serfaty \cite{AmSr}. The understanding of this limit has been done in collaboration with Sylvia Serfaty \cite{SerVaz}, and is related to work by Bertozzi et al. on aggregation models \cite{Bertozzi1, BLL}.

The more general equation $u_t=\nabla (u\, \nabla {\mathcal K} (u^{m-1}))$, $m>0$, is considered in \cite{BIK2013}, and self-similar solutions are found in explicit form, a very interesting fact.  A very natural extension of these equations is
\begin{equation}
u_t=\nabla (u^{m-1}\, \nabla {\mathcal K} (u^p)),\qquad
\end{equation}
where we may take for simplicity $u\ge 0$, $m>1$ and $p>0$. We are working in such generalizations, cf. \cite{StanTesoVazquez}  where the property of finite speed of propagation is examined.

\smallskip

The study of the fine asymptotic  behavior, i.e., obtaining rates of convergence is not easy. Progress on this issue is under way in collaboration with Carrillo and Huang \cite{CH2013}.

\smallskip

\noindent $\bullet$ {\bf Disclaimer.} Since our objective is limited, we are leaving completely outside of the presentation many related topics. Thus, the anomalous diffusion that is  encountered in many applications  is often described by means of {\sl fractional time operators}, using the so-called Riemann-Liouville derivatives.  There is a huge literature on the subject in the case of linear equations. We will not enter into such topic in this document, though some nonlinear models are quite appealing.

\section{The fractional porous medium equation}
\label{sec.fpme}

We now turn our attention to the nonlinear heat equation with fractional diffusion
\eqref{fpme}. Indeed, it is a whole family of  equations with exponents $s\in (0,1)$ and $m>0$. They can be seen as fractional-diffusion versions of the PME described above, \cite{Vapme}, \cite{JLVSmoothing}. The classical Heat Equation  is recovered in this model in the limit $s=1$ when $m=1$, the PME when $m>1$, the Fast Diffusion Equation when $m<1$.

\subsection{Some applied literature}\label{App.Motiv}

We gather here some updated information on the occurrence of the nonlinear fractional diffusion equation we propose and related models in the physical or probabilistic literature.

\smallskip

\noindent $\bullet$  Anomalous diffusion often takes a nonlinear form. To be more specific, there exist many phenomena  in nature where, as time goes on, a crossover is observed between different diffusion regimes. Tsallis et al. \cite{BGT2000, LMT2003}  discuss the following cases: (i) a mixture of the porous medium equation, which is connected with non-extensive statistical mechanics, with the normal diffusion equation; (ii) a mixture of the fractional time derivative and normal diffusion equations; (iii) a mixture of the fractional space derivative, which is related with L\'evy flights, and normal diffusion equations. In all three cases a crossover is obtained between anomalous and normal diffusions. This leads to models of nonlinear diffusion of porous medium or fast diffusion types with standard or  fractional Laplace operators, cf.  equation (4) of \cite{BGT2000}.

\smallskip

\noindent $\bullet$  There have been many studies of hydrodynamic limits of interacting particle systems with long-range dynamics, which lead to fractional diffusion equations of our type, mainly linear like in \cite{JKOlla}, but also nonlinear in the recent literature, cf. the works \cite{Jara1}, \cite{Jara2}. Thus, in the last reference, Jara and co-authors study the non-equilibrium functional central limit theorem for the position of a tagged particle in a mean-zero one-dimensional zero-range process. The asymptotic
behavior of the particle is described by a stochastic differential equation
governed by the solution of the following nonlinear hydrodynamic (PDE) equation, $\partial_t \rho = a^2 \partial^2_x \Phi(\rho)$. When $\Phi$ is a power we recover equation \eqref{fpme}.

\smallskip

\noindent $\bullet$ Equations like the last one (in several space dimensions) occur in boundary heat control, as already mentioned by Athanasopoulos and Caffarelli \cite{AC09}\,, where they refer to the model  formulated in the book by Duvaut and Lions \cite{DL1972}, and use the so-called Caffarelli-Silvestre extension \cite{Caffarelli-Silvestre}.

\subsection{Mathematical problem and general notions}
\label{sect-intro}

Let us present the main features and results in the theory we have developed.
To be specific, the theory of existence and
uniqueness as well the main properties are first studied by De Pablo,
Quir\'{o}s,  Rodr\'{\i}guez, and  V\'{a}zquez in \cite{pqrv1,
pqrv2} for the Cauchy problem
\begin{equation}  \label{eq:main}
\left\{
\begin{array}{ll}
\dfrac{\partial u}{\partial t} + (-\Delta)^{\sigma/2}
(|u|^{m-1}u)=0, & \qquad  x\in\ren,\; t>0,
\\ [4mm]
u(x,0) = f(x), & \qquad x\in\ren.%
\end{array}
\right.
\end{equation}
We have put $\sigma=2s$. The notation $|u|^{m-1}u$ is used to allow for solutions of two signs, but since it is a bit awkward we will often write  $u^m$ instead of $|u|^{m-1}u$ even for signed solutions when no confusion is feared; the same happens with the power $u^{1/m}$. We take initial data $f\in L^1(\ren)$, which is a standard assumption in diffusion problems on physical grounds. As for the exponents, we consider the whole fractional exponent range $0<\sigma<2$, and take porous medium exponent $m>0$.  As we have said, in the limit $\sigma\to 2$ we want to recover the standard Porous Medium Equation (PME) \ $u_t-\Delta (|u|^{m-1}u)=0.$

\smallskip

The two papers contain a rather complete analysis of the problem. A semigroup of weak energy solutions is constructed for every choice of $m$ and $\sigma$, both the $L^\infty$ smoothing effect  and $C^\alpha$ regularity work in most cases (with some restrictions if $m$ is near 0), and there is infinite propagation for all $m>0$ and $\sigma<1$.  The results can be viewed as a nonlinear interpolation between  the extreme cases $\sigma=2$: \ $u_t -\Delta(|u|^{m-1}u)=0$,
and $\sigma=0$ which turns out to be a simple ODE: \ $u_t+|u|^{m-1}u=0$.
It is to be noted that the critical exponent $m_c:= (N-\sigma)_+/N$   plays a role
in the qualitative theory: the properties of the semigroup are more familiar when $m>m_c$.
A similar exponent is well-known in the Fast Diffusion theory (putting $\sigma=2$). Note that
such exponent is not considered when $N=1$ and $\sigma\ge 1$.

\smallskip

\noindent{\bf Preliminary notions. } If $\psi$ and $\varphi$ belong to the Schwartz class, the definition \eqref{def-fourier} of the fractional Laplacian together with
Plancherel's theorem yield
 $$
\int_{\ren}(-\Delta)^{\sigma/2}\psi\,\varphi=\int_{\ren}|\xi|^\sigma
  \widehat\psi\,\widehat\varphi=\int_{\ren}|\xi|^{\sigma/2} \widehat
\psi|\xi|^{\sigma/2}\,\widehat
\varphi=\int_{\ren}(-\Delta)^{\sigma/4}\psi\,(-\Delta)^{\sigma/4}\varphi.
$$
Therefore, if we multiply  the equation in \eqref{eq:main} by a test
function $\varphi$ and integrate by parts, we obtain
\begin{equation}\label{weak-nonlocal}
\displaystyle \int_0^T\int_{\ren}u\dfrac{\partial
\varphi}{\partial
t}\,dxds-\int_0^T\int_{\ren}(-\Delta)^{\sigma/4}(|u|^{m-1}u)(-\Delta)^{\sigma/4}\varphi\,d
xds=0.
\end{equation}
This identity is the basis of our definition of a weak solution.
The integrals in \eqref{weak-nonlocal} make sense if $u$ and $u^m$
belong to suitable spaces. The right space for $u^m$ is the
fractional Sobolev space $\dot{H}^{\sigma/2}(\ren)$, defined
as the completion of $C_0^\infty(\ren)$ with the norm
$$
  \|\psi\|_{\dot{H}^{\sigma/2}}=\left(\int_{\ren}
|\xi|^\sigma|\widehat{\psi}|^2\,d\xi\right)^{1/2}
=\|(-\Delta)^{\sigma/4}\psi\|_{L^2},
$$
where $\widehat{\psi}$ denotes the Fourier transform of $\psi$. Note that $(-\Delta)^{\sigma/4}u^m\in L^2(\ren)$ if $u^m\in\dot{H}^{\sigma/2}(\ren)$.

\medskip

\noindent{\bf Definition}\label{def:weak.solution.nonlocla} A function $u$
is a {\it weak solution} to Problem \eqref{eq:main} if:
\begin{itemize}
\item $u\in
L^{1}(\ren\times (0,T))$ for all $T>0$, $u^m \in L^2_{\rm
loc}((0,\infty);\dot{H}^{\sigma/2}(\ren))$;

\item  identity  \eqref{weak-nonlocal}
holds for every $\varphi\in C_c^1(\ren\times(0,T))$;
\item  $u(\cdot,t)\in L^1(\ren)$ for all $t>0$,
$\lim\limits_{t\to0}u(\cdot,t)=f$ in $L^1(\ren)$.
\end{itemize}

A drawback of this definition is that there is no convenient formula
for the fractional Laplacian of a product or of a composition of
functions. Moreover, we take no advantage in using compactly
supported test functions since their fractional Laplacian loses this
property. To overcome these and other difficulties, we will use the fact that
our solution $u$ is the trace of the solution of a \emph{local}
problem obtained by extending $u^m$ to a half-space whose boundary
is our original space.

\smallskip

\noindent {\sc Extension Method.}  In  the particular case
$\sigma=1$ studied in \cite{pqrv1}, the problem is reformulated  by means of  the well-known representation of the
half-Laplacian in terms of the Dirichlet-Neumann operator. This allowed us to transform the nonlocal problem into a local one
(i.\,e., involving only derivatives and not integral operators). Of
course, this simplification pays a prize, namely, introducing an
extra space variable. The application of such an idea is not so
simple when $\sigma\ne 1$; it involves a number of difficulties that
we address in \cite{pqrv2}.  We have to use  the
characterization of the Laplacian of order $\sigma$,
$(-\Delta)^{\sigma/2}$, $0<\sigma<2$, described by
Caffarelli and Silvestre in their famous paper \cite{Caffarelli-Silvestre}, in terms of
the so-called $\sigma$-harmonic extension, which is the solution of
an elliptic problem with a degenerate or singular weight.

\smallskip

Let  us explain this extension in  some more detail. If $g=g(x)$ is
a smooth bounded function defined in $\ren$, its
$\sigma$-harmonic extension  to the upper half-space
$\mathbb{R}^{N+1}_+$, $v=~\E(g)$, is  the unique smooth bounded
solution $v=v(x,y)$ to
\begin{equation}
\left\{
\begin{array}{ll}
\nabla\cdot(y^{1-\sigma}\nabla v)=0,\qquad &x\in\ren,\, y>0,\\
v(x,0)=g(x),\qquad&x\in\ren.
\end{array}
\right. \label{sigma-extension}\end{equation} Then,
\begin{equation}
-\mu_{\sigma}\lim_{y\to0^+}y^{1-\sigma}\frac{\partial v}{\partial
y}=(-\Delta)^{\sigma/2} g(x), \label{fract-lapla}
\end{equation}
where  the precise constant, which does not depend on $N$, is
$\mu_{\sigma}=\frac{2^{\sigma-1}\Gamma(\sigma/2)}{\Gamma(1-\sigma/2)}$,
see \cite{Caffarelli-Silvestre}.  In \eqref{sigma-extension} the
operator  $\nabla$ acts in all $(x,y)$ variables, while in
\eqref{fract-lapla} $(-\Delta)^{\sigma/2}$ acts only on the
$x=(x_1,\cdots,x_N)$ variables.  In the sequel we denote
$$
  L_\sigma v\equiv \nabla\cdot(y^{1-\sigma}\nabla v),\qquad
  \dfrac{\partial v}{\partial y^\sigma}\equiv
  \mu_{\sigma}\lim_{y\to0^+}y^{1-\sigma}\frac{\partial v}{\partial
  y}.
$$

\noindent\emph{Notation.} In this section we will use the notation $\Omega=\ren\times(0,\infty)$ for the upper half-space, with points
$\overline x=(x,y)$, $x\in\ren$, $y>0$; its boundary, which is identified to the original
$\ren$  with variable $x$, will be named $\Gamma$.
Besides, we use the simplified notation $u^m$  for data of any sign,  instead of the
actual  ``odd power'' $|u|^{m-1}u$, and we will also use such a notation
when $m$ is replaced by $1/m$. The convention is not applied to any
other powers.

\smallskip

\noindent {\sc Extended problem. Weak solutions.} With the above in
mind, we rewrite problem \eqref{eq:main} for $w=u^m$ as a
quasi-stationary problem with a dynamical boundary condition
\begin{equation}
\left\{
\begin{array}{ll}
L_\sigma w=0\qquad &\mbox{for } \overline x\in\Omega,\, t>0,\\
\dfrac{\partial w}{\partial y^\sigma}-\dfrac{\partial
w^{1/m}}{\partial
t}=0\qquad&\mbox{for } x\in\Gamma,\, t>0,\\
w(x,0,0)=f^m(x)\qquad&\mbox{for } x\in\Gamma.
\end{array}
\right. \label{pp:local}
\end{equation}
This problem has  been considered  by Athanasopoulos and Caffarelli
\cite{AC09}, who prove that any bounded weak
solution is H\"{o}lder continuous if $m>1$.

\smallskip

To define a weak solution of this problem we multiply formally the
equation in \eqref{pp:local} by a test function $\varphi$ and
integrate by parts to obtain
\begin{equation}\label{weak-local}
\displaystyle \int_0^T\int_{\Gamma}u\dfrac{\partial
\varphi}{\partial
t}\,dxds-\mu_\sigma\int_0^T\int_{\Omega}y^{1-\sigma}\langle\nabla
w,\nabla \varphi\rangle\,d\overline xds=0,
\end{equation}
where $u=(\mathop{\rm Tr}(w))^{1/m}$ is the trace of $w$ on $\Gamma$ to the
power $1/m$. This holds on the condition that $\varphi$ vanishes for
$t=0$ and $t=T$, and also for large $|x|$ and $y$. We then introduce
the energy space $X^\sigma(\Omega)$, the completion of
$C_0^\infty(\Omega)$ with the norm
\begin{equation}
  \|v\|_{X^\sigma}=\left(\mu_\sigma\int_{\Omega} y^{1-\sigma}|\nabla
v|^2\,d\overline x\right)^{1/2}. \label{norma2}\end{equation} The
trace operator is well defined in this space, see below.

\smallskip

\noindent{\bf Definition}\label{def:weak.solution}
A pair of functions $(u,w)$ is a {\it weak solution} to Problem
\eqref{pp:local} if:
\begin{itemize}
\item $w\in L^2_{\rm loc}((0,\infty);X^\sigma(\Omega))$, $u=(\mathop{\rm Tr}(w))^{1/m}\in
L^{1}(\Gamma\times (0,T))$ for all $T>0$;

\item  Identity  \eqref{weak-local}
holds for every $\varphi\in C_0^1(\overline\Omega\times(0,T))$;
\item  $u(\cdot,t)\in L^1(\Gamma)$ for all $t>0$,
$\lim\limits_{t\to0}u(\cdot,t)=f$ in $L^1(\Gamma)$.
\end{itemize}

\smallskip

For brevity we will refer sometimes to the solution as only $u$, or
even only $w$, when no confusion arises, since it is clear how to
complete the pair from one of the components, $u=(\mathop{\rm Tr}(w))^{1/m}$,
$w=\E(u^m)$.

\smallskip

\noindent{\sc Equivalence of weak formulations.} The key point of
the above discussion is that the definitions of weak solution for
our original nonlocal  problem and for the extended local problem
are equivalent. The main ingredient of the proof is that
equation~\eqref{fract-lapla} holds in the sense of distributions for
any $g\in \dot H^{\sigma/2}(\Gamma)$.

\smallskip

\noindent {\bf Proposition} {\sl A function  $u$ is a weak solution to Problem \eqref{eq:main} if and only if $(u,\E(u^m))$ is a weak solution to Problem~\eqref{pp:local}.}

\medskip

\noindent{\sc Strong solutions.} Weak solutions satisfy
equation~\eqref{eq:main} in the sense of distributions. Hence, if
the left hand side is a function, the right hand side is also a
function and the equation holds almost everywhere. This fact allows
to prove uniqueness and several other important properties, and
hence motivates the following definition.

\medskip

\noindent{\bf Definition} We say that a weak solution $u$  to Problem \eqref{eq:main} is a strong solution
if $u\in C([0,\infty):\,L^1(\Gamma))$ \ as well as \ $\partial_tu$ and $(-\Delta)^{\sigma/2}
(|u|^{m-1}u) \in
L^1_{\rm loc}(\Gamma\times(0,\infty))$.

\subsection{Main results}

\noindent{\sc Existence.} We prove existence of a suitable
concept of (weak) solution for general $L^1$ initial data only in
the restricted range $m>m_c$, which includes as
a particular case the linear fractional heat equation, case $m=1$.
If $0<m\le m_c$ (which implies that $0<\sigma<1$ if $N=1$, recall that $m_c=(N-\sigma)_+/N$) we need to slightly restrict the data to obtain weak solutions.

\begin{theorem}\label{th:existence}
If either $f\in L^1(\ren)$ and $m>m_c$, or $f\in
L^1(\ren)\cap L^p(\ren)$  with  $p>p_*(m)=(1-m)n/\sigma$  and
 $ 0<m\le m_c $, there exists a  weak solution to the Cauchy problem for the FPME.
\end{theorem}

\noindent {\sc Uniqueness.} We first prove uniqueness of solutions in the range $m\ge m_{c}$. If $0<m<m_{c}$, we need to use the concept of strong solution, a concept that is standard in the abstract theory of evolution equations. This is no restriction in view of the next results proved in \cite{pqrv2}.

\begin{theorem}\label{th:strong}
The solution given by Theorem~{\rm \ref{th:existence}} is a strong
solution.
\end{theorem}
We state the uniqueness result in its simplest version.

\begin{theorem}\label{th:uniqueness} For every $f$ and $m>0$  there exists at most one strong solution to Problem \eqref{eq:main}.
\end{theorem}

\noindent {\sc Qualitative behaviour.} The solutions to Problem
\eqref{eq:main} have some nice properties that are summarized here.

\begin{theorem}\label{th:properties}
Assume $f,f_1,f_2$ satisfy the hypotheses of Theorem~{\rm
\ref{th:existence}}, and let  $u,u_1,u_2$ be the corresponding
strong solutions to Problem~\eqref{eq:main}.

\emph{(i)} If $m\ge m_c$, the mass $\int_{\ren}u(x,t)\,dx$ is
conserved.

\emph{(ii)} If $0<m< m_c$, then $u(\cdot,t)$ vanishes identically in
a finite time $T(u_0)$.

\emph{(iii)} A smoothing effect holds in the form:
\begin{equation}
\|u(\cdot,t)\|_{L^\infty(\ren)}\le C\,t^{-\alpha_p}
\|f\|_{L^p(\ren)}^{\delta_p} \label{eq:L-inf-p}
\end{equation}
with $\alpha_p=(m-1+{\sigma}p /N)^{-1}$, $\delta_p=\sigma
p\alpha_p/N$, and \ $C=C(m,p,N,\sigma)$.  This holds for all $p\ge 1$ if
$m>m_c$, and only for $p>p_*(m)$ if $0<m\le m_c$.

\emph{(iv)} Any $L^p$-norm  of the solution, $1\le p\le \infty$, is
nonincreasing in time.

\emph{(v)} There is an $L^1$-order-contraction property,
$$
\int_{\ren}(u_1-u_2)_+(x,t)\,dx\le
\int_{\ren}(u_1-u_2)_+(x,0)\,dx.
$$

\emph{(vi)} If $f\ge0$ the solution is positive  for all $x$ and all
positive times if $m\ge m_c$ (resp. for all $x$ and all $0<t<T$ if
it vanishes in finite time $T(u_0)$ when $0<m<m_c$).

\emph{(vii)} If either  $m\ge 1$ or $f\ge0$, then $u\in
C^\alpha(\ren\times(0,\infty))$ for some $0<\alpha<1$.
\end{theorem}

In the linear case $m=1$ the above properties: conservation of mass,
the smoothing effect with a precise decay rate, positivity and
regularity, can be derived directly from the representation formula
\eqref{lineal} and the properties of the kernel $K_\sigma$.

\smallskip

\noindent {\sc Continuous dependence.} We also show that the solution (i.e., the semigroup) depends continuously  on the initial data and on both parameters $m$ and $\sigma$, in
particular in the nontrivial limit $\sigma\to 2$, that allows to recover the
standard PME, $\partial_t u-\Delta |u|^{m-1}u=0$, or the other end $\sigma\to
0$, for which we get the ODE : $\partial_t u+ |u|^{m-1}u=0$. Continuity will be
true in general only in $L^1_{\rm loc}$, unless we stay in the region of
parameters where mass is conserved.

\begin{theorem}\label{th:contdep}
The strong solutions depend continuously  in the norm of the space
$C([0,T]:\,L^1_{\rm loc}(\ren))$ on the parameters $m$,
$\sigma$, and the initial data $f$. If moreover $m\ge m_c$ and
$0<\sigma\le2$, convergence also holds in
$C([0,T]:\,L^1(\ren))$.
\end{theorem}

We will extend this theory in different ways. See in particular Section \ref{sec.gnl}.

 \section{Self-similar solutions and their role}\label{sec.bar}

Several types of special solutions play a role in the theory of the FPME: fundamental solutions, very singular solutions, extinction solutions,...

 \subsection {\bf Fundamental solutions}
 Let us consider the Cauchy problem for the FPME taking as initial data a Dirac delta,
\begin{equation}\label{eq.id}
u(x,0) = M\delta(x)\, \qquad M>0.
\end{equation}
Solutions with such data are called {\sl fundamental solutions} in the linear theory, and we will keep that name though their relevance is different in the nonlinear context.
We use the concept of  continuous and nonnegative weak solution introduced in \cite{pqrv1}, \cite{pqrv2}, for which there is a well-developed theory when the datum is a  function in $L^1(\ren)$ as we have outlined. We recall that the question under discussion has a well-known answer for the standard Porous Medium Equation, i.e., the limit case  $s=1$, in the form of the usual Barenblatt solutions; they were actually discovered in the 1950's by Zeldovich-Kompanyeets \cite{ZK50} and Barenblatt \cite{Bar52}, cf. their use in applications in \cite{BarBook} and in theory of the equation in \cite{Vapme}.

\smallskip

Let us describe the results of our paper \cite{VazBar2012}, devoted to prove  existence, uniqueness and main properties of the fundamental solutions for the FPME. The exponent $m$ varies in principle in the range $m>1$, but the methods extend to the linear case $m=1$ and even to the fast diffusion range (FDE)  $m<1$ on the condition that $m> m_c=\max\{(N-2s)/N,0\}$. Such a type of restriction on $m$ from below carries over from the PME-FDE theory, \cite{JLVSmoothing}. The  solutions we find are self-similar functions of the form \
\begin{equation}\label{sss}
u^*(x,t)= t^{-\alpha} F(|x|\,t^{-\beta})
\end{equation}
with suitable exponents $\alpha$ and $\beta$ and profile $F\ge 0$. Here is the main result.

\begin{theorem} \label{thm.Bs} For every choice of parameters $s\in(0,1)$ and $m>m_c$ where $m_c=\max\{(N-2s)/N,0\}$, and for every $M>0$, equation \eqref{fpme} admits a unique fundamental solution $u^*_M(x,t)$; it is a nonnegative and continuous weak solution for \ $t>0$ and takes the initial data \eqref{eq.id} as a trace in the sense of Radon measures. It has the self-similar form \eqref{sss} for suitable $\alpha$ and $\beta$ that can be calculated in terms of $N$ and $s$ in a dimensional way, precisely
\begin{equation}\label{scale.expo}
\alpha=\frac{N}{N(m-1)+2s}, \qquad \beta=\frac{1}{N(m-1)+2s}\,.
\end{equation}
The profile function  $F_M(r)$, $r\ge 0$, is a  bounded and continuous function, it is positive everywhere, it is monotone and  it goes to zero at infinity.
\end{theorem}

Let us point out that the Brownian scaling ($\beta=1/2$ and $\alpha=N/2$) is recovered in the fractional case under the condition $N(m-1)=2(1-s)$, for instance for $m=(N+1)/N$ in the most common case $s=1/2$. Moreover, the profile $F_M$ can be obtained by a simple rescaling of $F_1$. $F$ is actually a smooth function by the regularity results of \cite{vpqr}.

\smallskip

 The initial data are taken in the weak sense of measures
\begin{equation}
\lim_{t\to 0}\int_{\ren} u(x,t)\phi(x)\,dx= M\phi(0)
\end{equation}
for all $\phi\in C_b(\ren)$, the space of continuous and bounded functions in $\ren$.
We will call these self-similar solutions  of Problem \eqref{fpme}-\eqref{eq.id} with given $M>0$ the {\sl Barenblatt solutions} of the fractional diffusion model by analogy with the PME and other prominent studied cases. The form of the exponents explains the already mentioned restriction on $m$ from below.

\smallskip

Idea of the proofs: existence of the Barenblatt solutions is done by approximation of the Cauchy problem with smooth initial data and passing to the limit, using suitable a priori estimates.  Uniqueness needs to go over to the potential equation as follows: we take the convolution of $u(x,t)$ with the Riesz kernel and define:
$$
U(x,t)= C_{n,s}\int \frac{u(y,t)}{|x-y|^{N-2s}}\,dy \,.
$$
Then $(-\Delta )^s U= u$, and using the equation we get equation (PE):
\begin{equation*}\label{eq.pot1}
U_t=((-\Delta )^{-s} u)_t=-(-\Delta )^{-s}(-\Delta )^{s}u^m=-u^m.
\end{equation*}
in other words, $U_t+((-\Delta )^s U)^m=0$. Though awkward-looking, this `dual equation' admits a good uniqueness theorem, as shown in \cite{VazBar2012}.

\smallskip

The properties of the profile function $F(r)$ are quite important, in particular the asymptotic behavior. This is why in \cite{VazBar2012} we derive the elliptic equation that it satisfies:
\begin{equation}\label{sss.form}
(-\Delta )^{s}F^m =\alpha F +  \beta y\cdot \nabla F=\beta \nabla\cdot (yF)
\end{equation}
Moreover, putting $ s'=1-s$ we have
$$
\nabla (-\Delta )^{-s'} F^m=-\beta \, y\, F\,,
$$
which in radial coordinates gives
\begin{equation}
L_{s'} F^m(r)=\beta \int_r^\infty rF(r)dr\,,
\end{equation}
where $L_{s'}$ the radial expression of operator $(-\Delta)^{-{s'}}$.
Using this 1D integral equation, the following characterization for the behavior of the fundamental solution is obtained. The critical exponents are $m_c=(N-2s)/N$ and  $m_1=N/(N+2s)$, so that $m_c<m_1<1$.

\begin{theorem} For every $m>m_1$ we have the asymptotic estimate
\begin{equation}
C_1\,M^{\mu}\le F_M(r)\,r^{N+2s} \le C_2\,M^{\mu},
\end{equation}
where $M=\int F(x)\,dx$, $C_i=  C_i(m,N,s)>0$, and $\mu=(m-m_1)(N+2s)\beta$. On the other hand, for $m_c<m<m_1$, there is a constant $C_\infty(m,N,s)$ such that
\begin{equation}
 F_M(r)\,r^{2s/(1-m)}=C_\infty.
\end{equation}
The case $m=m_1$ is borderline and has a logarithmic correction.
\end{theorem}

\smallskip

\noindent {\bf Solutions with explicit form.} As this survey is written, Y. Huang reports  \cite{Huang2013} the explicit expression of the Barenblatt solution for the special value of  $m$, $m_{ex} = (N+2-2s)/(N+2s)$. The profile is given by
\begin{equation}
F_M(y) =\lambda\,(R^2 + |y|^2)^{-(N+2s)/2}
\end{equation}
where the two constants $\lambda$ and $R$ are determined by the total mass $M$ of the solution  and the parameter $\beta$. Note that for $s=1/2$ we have $m_{ex}=1$, and the solution is the one mentioned in the introduction for the linear case, \eqref{kernel.lin}. We always have $m_{ex}>m_1$;  the range of $m$ that is covered for $0<s<1$ is $N/(N+2)<m_{ex}<(N+2)/N$.

\smallskip

\subsection{\bf Asymptotic behavior} As a main application of this construction, we  prove in \cite{VazBar2012} that the asymptotic behaviour as $t\to\infty$ of the class of nonnegative weak solutions of \eqref{fpme} with finite mass  (i.\,e., $\int_{\ren} u(x,t)\,dx<\infty$) is given in first approximation by the family $u_M^*(x,t)$, i.\,e., the Barenblatt solutions are the attractors in that class of solutions.

\begin{theorem}\label{thm.exlimit} Let $u$ be a nonnegative solution of the FPME with  initial data $u_0=\mu \in {\mathcal M}_+(\ren)$, and $m>m_c$. Let $M=\mu(\ren)$ and let $u^* _M$ be the self-similar Barenblatt solution with mass $M$. Then  as $t\to\infty$ the solutions $u(x,t)$ and $u^*_M(x,t)$ are increasingly similar, and more precisely we have
\begin{equation}\label{conver.express.1}
\lim_{t\to\infty} \|u(\cdot,t)-u^*_M(\cdot,t)\|_1=0\,,
\end{equation}
and also
\begin{equation}\label{conver.express}
\lim_{t\to\infty} t^{\alpha}\,|u(x,t)-u^*_M(x,t)|=0\,, \qquad \alpha=N/(N(m-1)+2s)\,,
\end{equation}
 uniformly in $x\in \ren$. It follows that for every $p\in (1,\infty)$ we have
 \begin{equation}\label{conver.express.p}
\lim_{t\to\infty}  t^{(p-1)\alpha/p}\|u(\cdot,t)-u^*_M(\cdot,t)\|_p=0\,.
\end{equation}
\end{theorem}

This result is a clear example of the important role that the fundamental solutions and their properties can play in the applications.

\subsection{\bf Very singular solutions}  Only in the range $m_c<m<m_1=N/(N+2s)$ we can pass to the limit  $M\to\infty$  and obtain a special solution that has a fixed isolated singularity at $x=0$ that has separate-variables form and we call very singular solution (VSS) by analogy  with the standard Fast Diffusion Equation.  The result represents a marked difference with the case $s=1$ where VSS exist in the larger range  $m_c<m<1$. It is another manifestation of the long-range interactions of the fractional Laplacian, that avoids some of the purely local estimates of the standard FDE with the classical Laplacian operator; indeed, extrapolation of the standard estimates would justify the existence of a VSS in cases where there is none.

\smallskip

Among other interesting qualitative properties of the equation, in \cite{VazBar2012} we prove an Aleksandrov reflection principle. This is related to symmetrization principles that we will discuss below.

\subsection{\bf Other solutions}

Self-similar solutions of the second kind, where the similarity exponents $\alpha$ and $\beta$ are only related by the compatibility condition $(m-1)\alpha +2s\beta=1$, but are otherwise free not to verify the constant-mass condition $\alpha= N\beta$ are constructed in the work with Volzone \cite{VazVol2}.

\smallskip

The question of extinction of mass in finite time that was introduced in \cite{pqrv2} and happens for $m<m_c<1$ is also best understood in terms of explicit solutions, which were constructed in \cite{VazVol2} and have the form
\begin{equation}\label{ext.solut}
U(x,t)^{1-m}=C\,\frac{T-t}{|x|^{2s}}.
\end{equation}

 \section{Estimates, better existence and positivity}\label{sec.mb}

In collaboration with M. Bonforte  we have studied the existence of quantitative a  priori estimates of a local type for  solutions of the FPME and we have then derived  consequences for the theory. Thus, in \cite{BV2012} we have dealt with the properties of nonnegative solutions of the Cauchy problem. Such kind of estimates were obtained for the standard PME by Aronson-Caffarelli \cite{ArCaff} and by the present authors for the standard FDE \cite{BV, BV-ADV}\,. The same local tools are not efficient for the present fractional model due to the nonlocal character of the diffusion operator, but then estimates occur in weighted spaces.

\subsection{Weighted $L^1$ estimates in the fast diffusion range}

We will concentrate on the case $m<1$. The use of suitable weight functions allows to prove crucial $L^1$-weighted estimates that enter substantially into the derivation of the main results. The results take different forms according to the value of the exponent $m$, a fact that is to be expected since it happens for standard FDE (i.e., the limit case $s=1$). When $s<1$ the equation is nonlocal, therefore we cannot expect purely local estimates to hold. Indeed, we will obtain estimates in weighted spaces if the weight satisfies certain decay conditions at infinity. We present first a technical lemma that shows the difference with the standard Laplacian.

\begin{lemma}\label{Lem.phi}
Let $\varphi\in C^2(\ren)$ and positive real function that is radially symmetric and decreasing in $|x|\ge 1$. Assume also that $\varphi(x)\le |x|^{-\alpha}$ and that $|D^2\varphi(x)| \le c_0 |x|^{-\alpha-2}$\,, for some positive constant $\alpha$ and for $|x|$ large enough. Then,  for all $|x|\ge |x_0|>>1$ we have
$$
|(-\Delta)^s\varphi(x)|\le \dfrac{c}{|x|^{\gamma}}\,,
$$
with $\gamma=\alpha+2s$ for  $\alpha<N$ and $\gamma=N+2s$ for  $\alpha>N$; the constant $c>0$ depends only on $\alpha,s,N$ and $\|\varphi\|_{C^2(\ren)}$. If  $\alpha>N$ the reverse estimate holds  from below  if $\varphi\ge0$: \ $|(-\Delta)^s\varphi(x)|\ge c_1 |x|^{-(N+2s)}$  for all $|x|\ge |x_0|>>1$\,.
\end{lemma}

A suitable particular choice is the function $\varphi$ defined for $\alpha>0$ as $\varphi(x)=1$ for $|x|\le 1$ and $\varphi(x)= \left(1+(|x|^2-1)^4\right)^{-\alpha/8}$ if $|x|\ge 1\,.$ We will use the following notations: $m_c=(N-2s)/N$, $m_1=N/(N+2s)$, $p_c=N(1-m)/2$,

\smallskip

These are the weighted $L^1$ estimates obtained in \cite{BV2012}.

\begin{theorem}\label{prop.HP.s}
Let $u\ge v$ be two ordered solutions to  the FPME with $0<m<1$. Let $\varphi_R(x)=\varphi(x/R)$ where $R>0$ and $\varphi$ is as in the previous lemma with $0\le \varphi(x)\le |x|^{-\alpha}$ for $|x|>>1$ and
$$
N-\frac{2s}{1-m}<\alpha< N+\frac{2s}{m}\,.
$$
Then, for all $0\le \tau,t <\infty$ we have
\begin{equation}\label{HP.s}
\begin{array}{l}
\left(\dint_{\ren}\big(u(t,x)- v(t,x)\big)\varphi_R\,dx\right)^{1-m}\le \\
\left(\dint_{\ren}\big(u(\tau,x)- v(\tau,x)\big)\varphi_R\,dx\right)^{1-m}
+ \dfrac{C_1 \,|t-\tau|}{R^{2s-N(1-m)}}
\end{array}
\end{equation}
with $C_1(\alpha,m,N)>0$.
\end{theorem}

It is remarkable that the  estimate holds for (very) weak solutions, even changing sign solutions.  Also, it is worth pointing out that the estimate holds both for $\tau<t$ and for $\tau>t$.  The estimate implies the conservation of mass when $m_c<m<1$, by letting $R\to \infty$. On the other hand,  when $0<m<m_c$ solutions corresponding to $u_0\in L^1(\ren)\cap L^p(\ren)$ with $p\ge N(1-m)/2s$\,, extinguish in finite time $T(u_0)>0$\,, (see e.g. \cite{pqrv2}); the above estimates provide a lower bound for the extinction time in such a case, just by letting $\tau=T$ and $t=0$ in the above estimates:
\begin{equation}\label{HP.s.T}
\frac{1}{C_1\,R^{N(1-m)-2s}}\left(\int_{\ren}u_0\,\varphi_R\,dx\right)^{1-m}\le T(u_0)\,.
\end{equation}
Moreover, if the initial datum $u_0$ is such that the limit as $R\to+\infty$ of the right-hand side diverges to $+\infty$, then the corresponding solution $u(t,x)$ exists (and is positive) globally in time.

\subsection{Existence of solutions in weighted $L^1$-spaces}
\label{sect.exist.large}

As a first consequence of these estimates, we can extend the $L^1$  existence theory of \cite{pqrv1, pqrv2} to an existence result for very weak solutions with non-integrable data in some weighted $L^1_{\varphi}$ space. In particular, bounded initial data or data with slow growth at infinity are allowed.

\begin{theorem}\label{exist.large}
Let $0<m<1$ and let $u_0\in L^1(\ren, \varphi\,dx)$, where $\varphi$ is as in Theorem $\ref{prop.HP.s}$ with decay at infinity $|x|^{-\alpha}$, $d-[2s/(1-m)]<\alpha<d+(2s/m)$. Then there exists a  very weak solution  $u(t,\cdot)\in L^1(\ren, \varphi\,dx)$ to the FPME in $[0,T]\times \ren$, in the sense that
\[
\int_0^T\int_{\ren}u(t,x)\psi_t(t,x)\,dx\dt
=\int_0^T\int_{\ren}u^m(t,x)(-\Delta)^s\psi(t,x)\,dx\dt
\]
for all $\psi\in C_c^\infty([0,T]\times\ren)\,.$ This solution is continuous in the weighted space, $u\in C([0,T]:L^1(\ren, \varphi\,dx))$\,.
\end{theorem}

Idea of the proof: We take $0\le u_{0,n}\in L^1(\ren)\cap L^\infty(\ren)$ be a non-decreasing sequence of initial data $u_{0,n-1}\le u_{0,n}$, converging monotonically to $u_0\in L^1(\ren, \varphi\dx)$\,, i.\,e., such that $\int_{\ren}(u_0- u_{n,0})\varphi\dx \to 0$ as $n\to \infty$. Consider the unique solutions $u_n(t,x)$ of the equation with initial data $u_{0,n}$. By the comparison results of \cite{pqrv2} we know that they form a monotone sequence. The weighted estimates \eqref{HP.s} show that the sequence is bounded in $L^1(\ren, \varphi\dx)$ uniformly in $t\in[0,T]$\,. By the monotone convergence theorem in $L^1(\ren, \varphi\dx)$, we know that the solutions $u_n(t,x)$ converge monotonically as $n\to \infty$ to a function $u(t,x)\in L^\infty ((0,T): L^1(\ren, \varphi\,dx))$. We then pass to the limit $n\to\infty$.\qed

\smallskip

\noindent\textbf{Remark. }The solutions constructed above only need to be integrable with respect to the weight $\varphi$, which has a tail of order less than $d+2s/m$. Therefore, we have proved existence of solutions corresponding to initial data $u_0$ that can grow at infinity as $|x|^{(2s/m)-\varepsilon}$ for any $\varepsilon >0$\,. Note that for the linear case $m=1$ this exponent is optimal in view of the representation of solutions in terms of the fundamental solution, but this does not seem to be the case for $m<1$.

\begin{theorem}[Uniqueness]
The solution constructed in Theorem $\ref{exist.large}$ by approximation from below is unique. We call it the minimal solution. In this class of solutions the standard comparison result holds, and also the estimates of Theorem $\ref{prop.HP.s}$.
\end{theorem}

\subsection{Lower bounds for fractional fast diffusion}

Section 3 of paper \cite{BV2012} studies the actual positivity of nonnegative solutions via quantitative lower estimates in the so called good fast diffusion range. We will use the
previous notation: $\beta:=1/[2s-N(1-m)]$\,, which is positive if $m>m_c$\,.

\smallskip

The following is the main quantitative estimate from below for positive solutions.

\begin{theorem}\label{thm.lower}
Let $R_0>0$, $m_c<m<1$ and let $0\le u_0\in L^1(\ren, \varphi\dx)$, where $\varphi$ is as in Theorem $\ref{prop.HP.s}$ with decay at infinity $O(|x|^{-\alpha})$, with \ $N-[2s/(1-m)]<\alpha<N+(2s/m)$. Let $u(t,\cdot)\in L^1(\ren, \varphi\dx)$ be a very weak solution to the FPME with initial datum $u_0$.  Then can calculate a time
\begin{equation}\label{t*}
t_*:=C_* \,R_0^{2s-N(1-m)}\,\|u_0\|_{ L^1(B_{R_0})}^{1-m}
\end{equation}
with quantified $C_*$ such that
\begin{equation}\label{low.1.thm}
\inf_{x\in B_{R_0/2}}u(t,x)\ge
K_1\,R_0^{-\frac{2s}{1-m}}\,t^{\frac{1}{1-m}}\quad \mbox{ if } \ 0\le t\le t_*\,,
\end{equation}
while
\begin{equation}
\inf_{x\in B_{R_0/2}}u(t,x)\ge K_2\, \|u_0\|_{L^1(B_{R_0})}^{2s\beta}
\,t^{-N\beta}\quad \mbox{ if } \ t\ge t_*\,.
\end{equation}
The positive constants $C_*,K_1,K_2$  depend only on $m,s$ and $N\ge 1$.
\end{theorem}

We point out that the diffusion in the equation is nonlocal, but the estimates are local. These estimates can be combined with the $L^\infty$ bounds from above of the papers \cite{pqrv1, pqrv2} to sandwich the solution from above and below. The paper continues to obtain sharp estimates on the behaviour as $|x|\to\infty$ (so-called tail behaviour).

\smallskip

\subsection{\bf Estimates in other ranges}
New local estimates are obtained in the paper either for $0<m<m_c$ or for $m\ge 1$. The question of Harnack inequalities is also discussed. We refrain from giving more details at this moment for lack of space, and refer to the paper \cite{BV2012} for further information.

\smallskip

An interesting result  is the existence and uniqueness of an initial trace for nonnegative solutions that we state in complete detail for $m<1$.

\begin{theorem}\label{thm.init.trace.m<1}
Let $0<m<1$ and let $u$ be a nonnegative weak solution of equation \eqref{fpme} in $(0,T]\times\ren$. Assume that $\|u(T)\|_{L^1(\ren)}<\infty$. Then there exists a unique nonnegative Radon measure $\mu$ as initial trace, that is
\begin{equation}\label{eq.trace1}
\int_{\ren}\psi\,\rd\mu=\lim_{t\to 0^+}\int_{\ren}u(t,x)\psi(x)\,\dx\,,\qquad\mbox{for all }\psi\in C_0(\ren)\,.
\end{equation}
Moreover, the initial trace $\mu$ satisfies the bound
\begin{equation}\label{intit.trace.bdd.lem}
\mu(B_{R}(x_0))\le \|u(T)\|_{L^1(\ren)} + C_1 R^{N(1-m)-2s}\,T\,.
\end{equation}
where $C_1=C_1(m,N,s)>0$ as in \eqref{HP.s}.
\end{theorem}

The result is indeed true for  $m>0$, including of course the linear equation for $m=1$.

 \section{Problems in bounded domains}\label{sec.dir}

The investigation of the questions of existence of solutions for the FPME posed in a bounded domain with suitable boundary conditions was initiated in the papers \cite{pqrv1, pqrv2}. The study of a priori estimates has been taken up in ongoing collaboration with
M. Bonforte. These estimates are quite different from the problem in the whole space in statements and methods.

\subsection{Fractional Laplacian operators on bounded domains}\label{ssect.Def.Fract.Lapl}

The first problem that we encounter is the possibility of various reasonable definitions of the concept of fractional Laplacian operator. This is in contrast with the situation in the whole Euclidean space $\re^N$, where there is a natural concept of fractional Laplacian that can be defined in several equivalent ways as we have mentioned in Subsection \ref{subs.frop}. When we consider  equation \eqref{fpme} posed on bounded domains, the definition via the Fourier transform does not apply, and different choices  appear as possible  definitions of the fractional Laplacian.

\smallskip

\noindent $\bullet$ On one hand, starting  from the classical Dirichlet Laplacian $\Delta_{\Omega}$ on the domain $\Omega$\,, the so-called spectral definition of the fractional power of $\Delta_{\Omega}$ uses the formula in terms of the semigroup associated to the Laplacian, namely
\begin{equation}\label{sLapl.Omega.Spectral}
\displaystyle(-\Delta_{\Omega})^{s}
g(x)=\frac1{\Gamma(-s)}\int_0^\infty
\left(e^{t\Delta_{\Omega}}g(x)-g(x)\right)\frac{dt}{t^{1+s}}.
\end{equation}
 We will call the operator defined in such a way  the \textit{spectral fractional Laplacian}, SFL for short, and denote  as $\A_1=(-\Delta_{\Omega})^s$. In this case, the initial and boundary conditions associated to the fractional diffusion equation \eqref{fpme} read
\begin{equation}\label{FPME.Dirichlet.conditions.Spectral}
\left\{
\begin{array}{lll}
u(t,x)=0\,,\; &\mbox{in }(0,\infty)\times\partial\Omega\,,\\
u(0,\cdot)=u_0\,,\; &\mbox{in }\Omega\,,
\end{array}
\right.
\end{equation}
Let us list some properties of the operator:  $\A_1=(-\Delta_{\Omega})^s$ is a self-adjoint operator on $ L^2(\Omega)$\,, with a discrete spectrum: its eigenvalues  are the family of $s$-power of the eigenvalues of the Dirichlet Laplacian: $(\lambda_j)^s>0$, $j=1,2,\ldots$; the corresponding normalized eigenfunctions \ $\Phi_{j}$ \   are exactly the same as the Dirichlet Laplacian, therefore they are as smooth as the boundary allows, namely when $\partial\Omega$ is $C^k$, then  $\Phi_j\in C^{\infty}(\Omega)\cap C^k(\overline{\Omega})$\,.
We can thus write
\begin{equation}\label{sLapl.Omega.Spectral.2}
\displaystyle(-\Delta_{\Omega})^{s}
g(x)=\sum_{j=1}^{\infty}(\lambda_j)^s\, \hat{g}_j\, \phi_j(x)
\end{equation}
with \ $
\hat{g}_j=\int_\Omega g(x)\Phi_j(x)\dx$\,, and $\|\Phi_j\|_{ L^2(\Omega)}=1\,.$

\medskip

The definition of the Fractional Laplacian via the Caffarelli-Silvestre extension \cite{Caffarelli-Silvestre} has been extended to the case of bounded domains by Cabr\'e and Tan \cite{Cabre-Tan}  by using as extended domain the cylinder ${\mathcal C}=(0,\infty)\times \Omega$ in $\re^{N+1}_+$, and by putting zero boundary conditions on the lateral boundary of that cylinder. This definition enables to understand the boundary conditions in an easy way. It is proved that this definition is equivalent to the SFL. See also \cite{capella-d-d-s}.

\smallskip

\noindent $\bullet$ On the other hand, we can define a fractional Laplacian operator by using the integral representation \eqref{def-riesz} in terms of hypersingular kernels and ``restrict'' the operator to functions that are zero outside $\Omega$: we will denote the operator defined in such a way as $\A_2=(-\Delta_{|\Omega})^s$\,, and call it the \textit{restricted fractional Laplacian}, RFL for short. In this case, the initial and boundary conditions associated to the fractional diffusion equation \eqref{fpme} read
\begin{equation}\label{FPME.Dirichlet.conditions.Restricted}
\left\{
\begin{array}{lll}
u(t,x)=0\; &\mbox{in }(0,\infty)\times (\re^N\setminus \Omega)\,,\\
u(0,\cdot)=u_0\; &\mbox{in }\Omega\,.
\end{array}
\right.
\end{equation}
See more in \cite{BVdir2013}. We also refer to \cite{SV1} for a  discussion and references about the differences between the Spectral and the Restricted fractional Laplacian. The authors of \cite{SV1} call the second type simply Fractional Laplacian, but we feel that the absence of descriptive name leads to confusion.

\subsection{\bf Main results on the Dirichlet problem}

In our paper  \cite{BVdir2013} we choose to work with the Dirichlet spectral laplacian, ${\mathcal L}_1$. For a quite general class of nonnegative weak solutions to the above problem, we derive in absolute upper estimates up to the boundary of the form
\begin{equation}\label{Intro.Abs.Bdds}
u(t,x) \le K_2\, \dfrac{\Phi_1(x)^{\frac{1}{m}}}{t^{\frac{1}{m-1}}}\,,\qquad\forall t>0\,,\;\forall x\in \Omega\,.
\end{equation}
In particular, we observe that the boundary behaviour is dictated by $\Phi_1$\,, the first positive eigenfunction of ${\A}_1$\,, which behaves like the distance to the boundary at least when the domain is smooth enough. We will also prove standard and weighted instantaneous smoothing effects of the type
\begin{equation}\label{Intro.smoothing}
\sup_{x\in \Omega}u(t,x)\le \dfrac{K_4}{t^{N\vartheta_{1,1}}}\left(\int_\Omega u(t,x)\Phi_1(x)\dx\right)^{2s\vartheta_{1,1}}
\le\dfrac{K_4}{t^{N\vartheta_{1,1}}}\left(\int_\Omega u_0\Phi_1(x)\dx\right)^{2s\vartheta_{1,1}}\,,
\end{equation}
where $\vartheta_{1,1}=1/(2s+(N+1)(m-1))$\,. This is sharper than \eqref{Intro.Abs.Bdds} only for small times.  As a consequence of the above upper estimates, we derive a number of useful weighted estimates  and we also obtain backward in time smoothing effect of the form
\begin{equation}\label{Intro.smoothing.back}
\| u(t)\|_{ L^\infty(\Omega)} \le \frac{K_4}{t^{(N+1)\vartheta_{1,1}}}\left(1\vee \frac{h}{t}\right)^{\frac{2s\vartheta_{1,1}}{m-1}}\|u(t+h)\|_{ L^1_{\Phi_1}(\Omega)}^{2s\vartheta_{1,1}}\,,\qquad\mbox{$\forall t,h>0$}\,,
\end{equation}
which is quite surprising and has not been observed before to our knowledge.

\smallskip

We then pass to the question of lower estimates, the main result being the estimate
\begin{equation}\label{Intro.lower}
u(t_0,x_0) \ge L_1\, \dfrac{\Phi_1(x_0)^{\frac{1}{m}}}{t^{\frac{1}{m-1}}}\qquad\mbox{$\forall t\ge t_* > 0\,,\;\forall x_0\in \Omega$}\,,
\end{equation}
where  the waiting time $t_*$  has the explicit form
\[
t_*= L_0\left(\int_{\Omega}u_0\Phi_1\dx\right)^{-(m-1)}\,.
\]
Then we observe that the above estimates combine into global Harnack inequalities in the following form:  for all $t\ge t_*$
\begin{equation}\label{Intro.GHP}
H_0\,\frac{\Phi_1(x_0)^{\frac{1}{m}}}{t^{\frac{1}{m-1}}} \le \,u(t,x_0)\le H_1\, \frac{\Phi_1(x_0)^{\frac{1}{m}}}{t^{\frac{1}{m-1}}}\,.
\end{equation}
We also provide as a corollary the local Harnack inequalities of elliptic type.

\smallskip

All the constants in the above results are universal, in the sense that may depend only on $N, m, s$ and $\Omega$\,, but not on $u$. They also have an almost explicit expression, usually given in the proof. Actually, we have tried to obtain quantitative versions of the estimates with indication of the dependence of the relevant constants. In some cases the estimates are absolute, in the sense that they are valid independently of the (norm of the) initial data.

\smallskip

\noindent {\bf New method.} It is worth mentioning that the proofs are based on new ideas with respect to what was used in the standard nonlinear diffusion or in the previous section for fractional diffusion in the whole space (both in our papers and in the references). Thus, we exploit the functional properties of the linear operator as much as possible. More precisely, we use the Green function of the fractional operator  even in the definition of solution, and we make use of estimates on its behaviour in the proofs. A key ingredient is thus the knowledge of good estimates for the Green function, that we discuss at length in the paper.

\smallskip

A careful inspection of the proofs shows that the presented method would allow to treat a quite wide class of linear operators, an issue that we shall discuss  in a forthcoming paper \cite{BV-FPMEBDpaper2}.  Here, we have written everything referring to the concrete case of the spectral fractional Laplacian in order to keep the exposition clear and to focus on the main ideas, but the arguments are devised in view of the wider applicability.

\smallskip

\subsection{\bf Existence theory}  The paper is complemented with a brief presentation of the theory of existence and uniqueness of the class of very weak solutions that we use. This complements the basic theory developed before in  \cite{pqrv1,pqrv2} and then extended in \cite{BV2012}. We have no space to enter this intriguing issue here.

\section{The work on general linearities}\label{sec.gnl}

In collaboration with Pablo, Quir\'os and Rodr\'{\i}guez we have examined the question of regularity of the nonnegative solutions of the FPME. Since we had proved continuity and positivity, the equation is not more degenerate, at least at the formal level, and the solutions should be smooth if we recall what happens in the classical heat equation and fast diffusion case when $s=1$. Actually, in \cite{vpqr} we have decided to study the question for the more general class of equations
\begin{equation}  \label{eq:mainphi}
\partial_tu+(-\Delta)^{\sigma/2}\varphi(u)=0 \,.
\end{equation}
To be specific, we deal with the Cauchy problem posed in $Q=\ren\times(0,\infty).$
The constitutive function $\varphi$ is assumed to be at least continuous and nondecreasing. Further conditions will be introduced as needed. This generality was motivated by previous work on the model of logarithmic porous medium type equation with fractional diffusion
$$
\partial_tu+(-\Delta)^{1/2}\log(1+u)=0,
$$
which is described in detail in \cite{pqrv3} and has interesting connections with drift-diffusion models.

\smallskip

The existence of a unique weak solution to the Cauchy problem for the FPME
has been fully investigated in~\cite{pqrv1, pqrv2} for the case where
$\varphi$ is a positive power. This theory is extended in \cite{vpqr} under suitable conditions on $\varphi$ that cover powers in particular.

\smallskip

Regarding regularity, for the linear fractional heat equation it follows from explicit representation with a kernel that solutions are $C^\infty$ smooth and bounded for every
$t>0$, $x\in\ren$, under the assumption that the initial data are integrable.
In the  nonlinear case such a representation is not available. Nevertheless, we will still be able to obtain that the solution is smooth if the equation is \lq\lq uniformly parabolic'',
$0<c\le \varphi'(u)\le C<\infty$. Our first result establishes that if the nonlinearity is smooth enough, compared to the order of the equation, $\max\{1,\sigma\}$, then bounded weak solutions are indeed {\sl classical solutions}.

\begin{theorem}\label{th:main} Let $u$ be a bounded weak solution to \eqref{eq:mainphi}, and assume $\varphi\in C^{1,\gamma}(\mathbb{R})$,
$0<\gamma<1$, and $\varphi'(s)>0$ for every
$s\in\mathbb{R}$. If $1+\gamma>\sigma$, then
$\partial_tu$ and $(-\Delta)^{\sigma/2}\varphi(u)$ are H\"older continuous functions and~\eqref{eq:mainphi} is satisfied everywhere.
\end{theorem}
The precise regularity of the solution is determined by the regularity of the nonlinearity $\varphi$; see the paper for the details. Notice that the condition $\varphi'>0$ together with the boundedness of $u$ implies that the equation is uniformly parabolic.

\smallskip

The idea of the proof is as follows: thanks to the results of Athanasopoulos and Caffarelli~\cite{AC09}, we  know that bounded weak solutions
are $C^\alpha$ regular for some $\alpha\in(0,1)$. In order to improve this regularity we write the  equation~\eqref{eq:mainphi} as a fractional linear heat equation with a source term. This term is in principle not very smooth, but it has some good properties. To be precise, given  $(x_0,t_0)\in Q$, we have
\begin{equation}\label{eq:nonlin.lin}
\partial_tu+(-\Delta)^{\sigma/2} u=(-\Delta)^{\sigma/2} f,
\end{equation}
where
\begin{equation*}
\label{eq:linear-f}
f(x,t):=u(x,t)-\frac{\varphi(u(x,t))}{\varphi'(u(x_0,t_0))},
\end{equation*}
after the time rescaling $t\to t/\varphi'(u(x_0,t_0))$. A very delicate study of the linear theory, with special attention to the properties of the kernel of the fractional Laplacian, allows to show that solutions to the linear equation~\eqref{eq:nonlin.lin} have the same regularity as $f$.  Next, using the nonlinearity we observe that  $f$ in the actual right-hand side is more regular than $u$ near $(x_0,t_0)$. We are thus in a situation that is somewhat similar to the one considered by Caffarelli and Vasseur in \cite{CaffVass},  where they deal with an equation, motivated by the study of geostrophic equations, of the form
\begin{equation*}
\label{eq:caffa-vasseur}
\partial_tu+(-\Delta)^{1/2} u=\mathop{\rm div}( {\bf v}u),
\end{equation*}
where $\bf v$ is a divergence free vector. Comparing with \eqref{eq:nonlin.lin}, we see two differences: in their case $\sigma=1$, and the source term is local. These two differences will significantly complicate our analysis.

\smallskip

\noindent {\bf Singular and degenerate equations.} The hypotheses made in
Theorem~\ref{th:main} excludes  all the powers
$\varphi(u)=|u|^{m-1}u$ for $m>0$, $m\neq1$, since they are degenerate ($m>1$)
or singular ($m<1$) at the level $u=0$. Nevertheless, a close look at our proof shows
that  we may in fact get  a \lq\lq local'' result, see the paper. Therefore,  we get for these
nonlinearities (and also for more general ones) a    regularity result in the positivity (negativity) set of the solution that implies that bounded weak solutions which are either positive or negative are actually classical.

\smallskip

\noindent {\bf Higher regularity.} If $\varphi$ is $C^\infty$  we prove that solutions are $C^\infty$. The result will be a consequence of the regularity already provided by Theorem~\ref{th:main} plus a special result for linear equations with variable coefficients. The case $\sigma<1$ is even more involved since we first have to raise the regularity in space exponent from $\sigma$ to 1.

\begin{theorem}\label{th:main2}
Let $u$ be a bounded weak solution to equation \eqref{eq:mainphi}. If $\varphi\in C^\infty(\mathbb{R})$, $\varphi'>0$ in $\mathbb{R}$,  then $u\in C^\infty(Q)$.
\end{theorem}

\noindent {\bf Comments and extensions.}
Work on extending such results to problems in bounded domains, and to more general operators is in project. Regarding related literature, let us  remark that Kiselev et al.~~\cite{KisNazVol} give a proof of $C^\infty$ regularity of a class of periodic solutions of geostrophic equations in 2D with $C^\infty$ data. Their methods are completely different to the ones used by us.
On the other hand, Cifani and Jakobsen propose in \cite{CJ11} an alternative $L^1$ theory dealing with a more general class of nonlocal porous medium equations, including strong degeneration and convection. The quantitative study of continuous dependence has been taken up recently in work of Alibaud et al. \cite{acj2011, acj2013}.

 \section{Estimates via symmetrization}\label{sec.symm}

 In two papers  \cite{VazVol, VazVol2} in collaboration with B. Volzone we have used the techniques of Schwarz and Steiner symmetrization to obtain a priori estimates, in many cases  with best constants, for the solutions of the FPME, or its more general version  $u_t+(-\Delta)^{s}A(u)=0$,  posed in the whole space or in a bounded domain. The main results concern the case of power nonlinearities, the equation is posed in the whole space and the exponent $m<1$. Elliptic results are obtained in \cite{VazVol} as a preliminary to the parabolic theory. In order to keep up with previous sections we will use the letter $\varphi$ for the nonlinearity instead of the $A$ in the paper.

\smallskip

\noindent {\sc Motivation.} Symmetrization is a very old geometrical idea that has become nowadays a popular tool of obtaining a priori estimates for the solutions of different partial differential  equations,  notably  those of elliptic and parabolic type. The application of  Schwarz  symmetrization to obtaining a priori estimates for elliptic problems is already described in \cite{Wein62}.  The  standard  elliptic result refers to the solutions of an equation of the form
$$
Lu=f,  \qquad Lu=-\sum_{i,j} \partial_i(a_{ij}\partial_j u)\,,
$$
posed in a bounded domain $\Omega\subseteq \ren$;  the coefficients $\{a_{ij}\}$ are assumed to be bounded, measurable and satisfy the usual ellipticity
condition; finally, we take zero Dirichlet boundary conditions on the boundary $\partial\Omega$. The  classical analysis introduced by Talenti \cite{Talenti1} leads to pointwise comparison between (i) the symmetrized version (more precisely the spherical decreasing rearrangement)  of the actual solution of the problem $u(x)$ and (ii) the radially symmetric solution $v(|x|)$ of some radially symmetric model problem which is posed in a ball with the same volume as $\Omega$.
Sharp a priori estimates for the solutions are then derived. Extensions of this method to more general problems or related equations have led to a
copious literature.

\medskip

\noindent {\sc Parabolic version.} This pointwise comparison fails for parabolic problems and the appropriate concept is comparison of concentrations, cf. Bandle \cite{Bandle} and V\'azquez \cite{Vsym82}. The latter considers the evolution problems of the form
\begin{equation}\label{evol.pbm}
\partial_t u=\Delta \varphi(u), \quad u(0)=u_0,
 \end{equation}
where $\varphi$ a monotone increasing real function and $u_0$ is a suitably given initial datum which is assumed to be integrable. For simplicity the problem was
posed
for $x\in \ren$,  but bounded open sets can be used as spatial domains.

\medskip

\noindent {\sc Fractional operators.} Symmetrization techniques were first applied to PDEs involving fractional Laplacian operators in the paper \cite{BV}, where the linear elliptic case is studied.
In our first paper  \cite{VazVol} we were able to improve on that progress and combine it with the  parabolic ideas of \cite{Vsym82} to establish the relevant
comparison theorems based on symmetrization for linear and  nonlinear parabolic equations.  To be specific, we deal with equations of the form \
\begin{equation}\label{nolin.parab}
\partial_t u +(-\Delta)^{s}\varphi(u)=f, \qquad 0<s<1\,.
\end{equation}
 Let us describe the results in \cite{VazVol} concerning the idea of {\sl concentration comparison} in the case of the solutions to the Cauchy problem
\begin{equation} \label{eq.1}
\left\{
\begin{array}
[c]{lll}
u_t+(-\Delta)^{\sigma/2}\varphi(u)=f  &  & x\in\mathbb{R}^{N}\,,t>0%
\\[6pt]
u(x,0)=u_{0}(x) &  & x\in\mathbb{R}^{N}\,.
\end{array}
\right.
\end{equation}
Here, the nonlinearity $\varphi(u)$ is a nonnegative function, smooth on $\mathbb{R}_{+}$,
with $\varphi(0)=0$ and $\varphi'(u)>0$ for all $u>0$ (extended anti-symmetrically in the general two-signed theory). Special attention is paid to  cases of the form $\varphi(u)=u^m$ with $m>0$.In \cite{VazVol} we have obtained that a concentration comparison for solutions to \eqref{eq.1} holds \emph{only} when the nonlinearity $\varphi$ is a \emph{concave} function, while for \emph{convex} $\varphi$ a remarkable example is
constructed, showing that a failure of concentration comparison occurs (see \cite{VazVol}).

\begin{theorem}\label{Main comparison}
Let $u$ be the nonnegative mild solution to problem \eqref{eq.1} with $\sigma\in(0,2)$,
initial data $u_0\in L^1(\ren)$, $u_0\ge 0$, right-hand side $f\in L^1(Q)$ where $Q=\mathbb{R}^{N}\times (0,\infty)$, $f\ge 0$, and nonlinearity $\varphi(u)$ given by a
concave function with
$\varphi(0)=0$
and
$\varphi'(u)>0$ for all $u>0$. Let $v$ be the solution of the symmetrized problem
\begin{equation} \label{eqcauchysymm.f}
\left\{
\begin{array}
[c]{lll}%
v_t+(-\Delta)^{\sigma/2}\varphi(v)=f^{\#}(|x|,t)  &  & x\in\mathbb{R}^{N}\,, \ t>0,%
\\[6pt]
v(x,0)=u_{0}^{\#}(x) &  & x\in\mathbb{R}^{N},
\end{array}
\right. %
\end{equation}
where $f^{\#}(|x|,t)$ means the spherical rearrangement of $f(x,t)$ w.r. to $x$ for fixed time $t>0$. Then,
for all $t>0$ we have
\begin{equation}
u^\#(|x|,t)\prec v(|x|,t).\label{conccompa}
\end{equation}
In particular, we have $\|u(\cdot,t)\|_p \le\|v(\cdot,t)\|_p$ for every $t>0$ and every $p\in [1,\infty]$.
\end{theorem}
Moreover, the following corollary justifies a reasonable consequence: if the data of problem \eqref{eq.1} are less concentrated that those of
the symmetrized problem, so are the corresponding solutions.
\begin{corollary}\label{corollarycomp}With the same assumptions of Theorem {\rm \ref{Main comparison}}, suppose that $u$ is the solution to problem
\eqref{eq.1} and
$v$ solves
\begin{equation}\label{eqcauchysymm.f1}
\left\{
\begin{array}
[c]{lll}%
v_t+(-\Delta)^{\sigma/2}\varphi(v)=\widetilde{f}(|x|,t)  &  & x\in\mathbb{R}^{N}\,, \ t>0,%
\\[6pt]
v(x,0)=\widetilde{u}_{0}(x) &  & x\in\mathbb{R}^{N},
\end{array}
\right. %
\end{equation}
where $\widetilde{f}\in L^{1}(Q)$, $\widetilde{u}_{0} \in L^{1}(\mathbb{R}^{N})$ are nonnegative, radially symmetric decreasing functions with respect to $|x|$. If
\[
u_{0}^{\#}(|x|)\prec\widetilde{u}_{0}(|x|),\quad f^{\#}(|x|,t)\prec\widetilde{f}(|x|,t)
\]
for almost all $t>0$, then  the conclusion $u^\#(|x|,t)\prec v(|x|,t)$  still holds.
\end{corollary}

\smallskip

\noindent {\bf A priori estimates  and best constants. } We  use the parabolic comparison results of \cite{VazVol} to obtain precise a priori estimates for the solutions of equation \eqref{nolin.parab}. One of these estimates is the so-called $L^1$ into $L^\infty$ smoothing effect in the following form:
\begin{equation}
 \|u(\cdot,t)\|_\infty\le C\,\|u_0\|_1^{2s\beta}\,t^{-\alpha}.
\end{equation}
The estimates were obtained in \cite{pqrv1, pqrv2} (and in the non-power case in \cite{vpqr}). In the paper \cite{VazVol2} we obtain the precise exponents and, what is important, the best constant $C$ in the decay inequality. The calculation of best constants in functional inequalities is a topic of continuing interest in the theory of PDEs, both in the elliptic and evolution settings. Classical references to the calculation of best constants by symmetrization methods are Aubin and Talenti's computation of the best constants in the Sobolev inequality in \cite{Aubin76, Talenti2}. Our calculation of a priori estimates with exact exponents and best
constants is closely related to  the sharp decay estimate for solutions of the porous medium/fast diffusion equation in \cite{Vsym82,JLVSmoothing}. When treating the linear case $\varphi(u)=u$, the estimates  are called ultra-contractivity, see the book \cite{Davies1} where the importance of best constants is stressed for the applications in Physics.
As a further application of the comparison techniques, optimal estimates with initial data in Marcinkiewicz spaces are obtained.

\smallskip

An important critical exponent appears repeatedly in the paper as a lower bound,
\begin{equation}
m_c:=(N-2s)/N\,.\label{supFFD}
\end{equation}
Since we are assuming $m>0$ and $0<s<1$, it does not appear in dimension $N=1$ if $s\ge 1/2$. Thus, we study the question of deciding the possible extinction of solutions in the range $m<m_c$ in terms of some norm of the initial data, and estimating  the
extinction time. First of all, we construct an explicit extinction solution of the fractional fast diffusion equation in this range of $m$'s, formula
\eqref{ext.solut}. Then, we obtain optimal estimates by using comparison based on symmetrization. In this direction we improve significantly the results  of  the previous papers \cite{BV2012} and \cite{pqrv2}, by obtaining optimal estimates on the extinction time for data in Marcinkiewicz spaces. Finally, all the  results are stable under the limit $s\to 1$, where the standard diffusion case is recovered. See \cite{pqrv2} for details on such limit.

 \section{ KPP propagation and fractional diffusion}\label{sec.kpp}

The problem goes back to the work of Kolmogorov, Petrovskii and Piskunov, see \cite{KPP}, that presents the most simple reaction-diffusion equation concerning the concentration $u$ of a single substance in one spatial dimension, $\partial_t u=D u_{xx} + f(u).$  The choice $f(u) = u(1-u)$ yields Fisher's equation \cite{Fisher} that was originally used to describe the spreading of biological populations. The celebrated  result says that the long-time behavior of any solution of $u(x,t)$, with suitable data $0\le u_0(x)\le 1$ that decay fast at infinity, resembles a traveling  wave with a definite speed.  In dimensions $N\geq 1$, the problem becomes
\begin{equation*}\label{classicalKPP}
u_t-\Delta u=u(1-u) \quad \text{in }(0,+\infty)\times \ren.
\end{equation*}
This case has been studied by Aronson and Weinberger in \cite{AronsonWeinberger}, where they prove the same result as the one-dimensional case. The result is formulated in terms of linear propagation of the level sets of the solution.
In the case of the more general model
\begin{equation}\label{PMEkpp}
u_t-\Delta u^m=u(1-u)
\end{equation}
the same result as before holds in the case of slow diffusion $m>1$ (V\'azquez and de Pablo \cite{dPJLVjde1991}). Departing from these results, King and McCabe  examined in \cite{KingMcCabe} the case of fast diffusion $m<1$ of equation \eqref{PMEkpp}. For $(N-2)_+/N<m<1$, the authors showed  that the problem does not admit traveling wave solutions by proving that level sets of the solutions of the initial-value problem with suitable initial data propagate exponentially fast in time.

\smallskip

On the other hand, and independently,  Cabr\'{e} and Roquejoffre (\cite{CabreRoquejoffre2}) studied the case of fractional linear diffusion $u_{t}(x,t) + (-\Delta)^s u(x,t)=f(u)$, where $(-\Delta)^s$ is the Fractional Laplacian operator with $s \in (0,1)$ and  concluded in the same vein that there is no traveling wave behavior as $t\to\infty$, and indeed the level sets propagate exponentially fast in time. This came as a surprise since their problem deals  with linear diffusion.

\smallskip

Motivated by these two examples of break of the asymptotic TW structure, we  studied in \cite{StanVazquezKPP} the case of a diffusion that is both fractional and nonlinear. More exactly, we consider the following reaction-diffusion problem
\begin{equation}\label{KPP}
  \left\{ \begin{array}{ll}
  u_{t}(x,t) + (-\Delta)^s u^m(x,t)=f(u) &\text{for } x \in \ren \text{ and }t>0, \\
  u(x,0)  =u_0(x) &\text{for } x \in \ren,
    \end{array}
    \right.
\end{equation}
We are interested in studying the propagation properties of nonnegative and bounded solutions of this problem in the spirit of the Fisher-KPP theory. We assume that the reaction term $f(u)$ satisfies
 \begin{equation*}\label{propf}
f\in C^1([0,1])\text{ is a concave function with } f(0)=f(1)=0 , \quad f'(1)<0<f'(0).
\end{equation*}
 For example we can take $f(u)=u(1-u).$ The initial datum $u_0(x) : \ren \rightarrow [0,1]$ and satisfies a growth condition of the
form
\begin{equation}\label{dataAssump}
0\le u_0(x) \leq C|x|^{-\lambda(N,s,m)}, \quad \forall x \in \ren,
\end{equation}
where the exponent $\lambda(N,s,m)$ is stated explicitly in the different ranges of the exponent $m$. In our work we establish the negative result about traveling wave behaviour,  more precisely, we prove that an exponential rate of propagation of level sets is true in all cases. Due to the nonlinearity, the solution of the diffusion problems involved in the proofs does not admit an integral representation as the case $m=1$. Instead, we use as an essential tool the behavior of the fundamental solution of the Fractional Porous Medium Equation, also called Barenblatt solution, recently  studied in \cite{VazBar2012}. This allows us to explain the propagation mechanism in simple terms: the exponential rate of propagation of the level sets of solutions (with initial data having a certain minimum decay for large $|x|$) is a consequence of the power-like decay behaviour of the fundamental solutions of the FPME.
To be precise, the decay rate of the tail of these solutions as $|x|\to\infty$ is the essential information we use to calculate the rates of expansion. This information is combined with more or less usual techniques of linearization and comparison with sub- and super-solutions.  We also need accurate lower estimates for positive solutions of this latter equation, and a further self-similar analysis for the linear diffusion problem. The delicate details are explained in \cite{StanVazquezKPP}.

\smallskip

\noindent {\bf Main results.} The existence of a unique mild solution of problem \eqref{KPP} follows by semigroup approach. The mild solution corresponding to an initial datum $u_0 \in
L^1(\ren),$ $0\leq u_0\leq1$ is in fact a positive, bounded, strong solution with $C^{1,\alpha}$ regularity. Let us introduce some notations. Once and for all, we put   $\beta=1/(N(m-1)+2s)$ \ and
\begin{equation}\label{sigma}
\sigma_1 =\frac{1-m}{2s}f'(0), \quad  \sigma_{2}=\frac{1}{N+2s}f'(0), \quad \sigma_3=\frac{1+2(m-1)\beta s}{N+2s}f'(0).
\end{equation}
The value $\sigma_1$ appears for $m_c<m<m_1$ and then $\sigma_1>\sigma_2$.
Notice also that  $\sigma_{2}<\sigma_3$ for $m>1$. Here is the precise statement of our main results for the solutions of the generalized KPP problem
\eqref{KPP}.

\begin{theorem}\label{mainThm1}
Let $N\geq 1$, $s\in (0,1)$, $f$ satisfying \eqref{propf} and $m_1<m \leq 1$. Let $u$ be a solution of \eqref{KPP}, where $0\leq u_0(\cdot)\leq 1$ is measurable,
$u_0 \neq 0$ and satisfies
\begin{equation}
0\le u_0(x) \leq C|x|^{-(N+2s)}, \quad \forall x \in \ren.
\end{equation}
 Then
\begin{enumerate}
\item if $\sigma>\sigma_{2}$, then $u(x,t) \rightarrow 0$ uniformly in $\{|x|\geq e^{\sigma t} \}$ as $t\rightarrow \infty$.
\item if $\sigma<\sigma_{2}$, then $u(x,t) \rightarrow 1$ uniformly in $\{|x|\leq e^{\sigma t} \}$ as $t\rightarrow \infty$.
\end{enumerate}
\end{theorem}

\begin{theorem}\label{mainThm2}
Let $N\geq 1$, $s\in (0,1)$, $f$ satisfying \eqref{propf} and $m_s<m<m_1$. Let $u$ be a solution of \eqref{KPP}, where $0\leq u_0(\cdot)\leq 1$ is measurable,
$u_0 \neq 0$ and satisfies
\begin{equation}
0\leq u_0(x) \leq C |x|^{-2s/(1-m)},  \quad \forall x \in \ren.
\end{equation}
  Then
\begin{enumerate}
\item if $\sigma>\sigma_{1}$, then $u(x,t) \rightarrow 0$ uniformly in $\{|x|\geq e^{\sigma t} \}$ as $t\rightarrow \infty$.
\item if $\sigma<\sigma_{1}$, then $u(x,t) \rightarrow 1$ uniformly in $\{|x|\leq e^{\sigma t} \}$ as $t\rightarrow \infty$.
\end{enumerate}
\end{theorem}

\begin{theorem}\label{mainThm3}
Let $N\geq 1$, $s\in (0,1)$, $f$ satisfying \eqref{propf} and $m>1$. Let $u$ be a solution of \eqref{KPP}, where $0\leq u_0(\cdot)\leq 1$ is measurable, $u_0 \neq
0$ and satisfies $$0\le u_0(x) \leq C|x|^{-(N+2s)}, \quad \forall x \in \ren.$$ Then
\begin{enumerate}
\item if $\sigma>\sigma_{3}$, then $u(x,t) \rightarrow 0$ uniformly in $\{|x|\geq e^{\sigma t} \}$ as $t\rightarrow \infty$.
\item if $\sigma<\sigma_{2}$, then $u(x,t) \rightarrow 1$ uniformly in $\{|x|\leq e^{\sigma t} \}$ as $t\rightarrow \infty$.
\end{enumerate}
\end{theorem}

Our main conclusion is  that  exponential propagation is shown to be the common occurrence, and the existence of traveling wave behavior is reduced to the classical KPP cases mentioned at the beginning of this discussion

\smallskip

Equations similar to \eqref{KPP} present interest for various researcher groups, like the groups led by X. Cabr\'e, J.M. Roquejoffre, H. Berestycki, F. Hamel, and P. Felmer.

 \section{Current work and comments}

 A number of related models, issues and perspectives on elliptic and parabolic equations involving fractional Laplacians  and more general integral operators is under way.

\smallskip

\noindent $\bullet $ {\bf Numerics.} The numerical analysis of  the solutions of the FPME was started by Teso \cite{Teso} where the case $s=1/2$ was studied. The extension method that we use to implement the fractional Laplacian has a number of specific difficulties for $s\ne 1/2$ that we have addressed in a subsequent paper  \cite{TesoVaz}. One of the main differences of our work is that we do not directly deal  with the integral formulation of the fractional  Laplacian;  instead of this, we pass through the  Caffarelli-Silvestre extension, mentioned above.

\smallskip

Previous works  dealing with the  numerical analysis of nonlocal equations of this type are due to Cifani, Jakobsen, and Karlsen in \cite{cjakobk}. In particular, they formulate some convergent numerical methods for entropy and viscosity solutions. The numerical analysis of the elliptic PDE $(-\Delta)^s u =f$ in a bounded domain with zero boundary data via the extension method has been recently studied by Nochetto and collaborators using  finite elements, \cite{NOS13}.

\

\noindent $\bullet$ {\bf Disclaimer.} We have not covered the extensive work on stationary states, i.e., elliptic equations of fractional type. Neither did we go into the probabilistic approach that comes from long ago and has been very actively pursued up to these days.
Evolution equations involving the fractional Laplacian appear in many applied models and some of this directions are now actively pursued, like  quasi-geostrophic flows, a topic that we have touched in the paper \cite{pqrv3}, but has not been addressed here. Another interesting direction concerns geometric flows with fractional operators. This list does not aim at being representative. Finally, we have not discussed the nonlocal diffusion evolution model treated in the monograph \cite{AMRT} which uses integral operators with kernels that do not resemble the fractional Laplacians and lead to a different theory.


\

{\sc Acknowledgments. } Work partially supported by Spanish Project MTM2011-24696. The author is very grateful to his main collaborators in the recent part of this effort: Matteo Bonforte, Arturo de Pablo, Fernando Quir\'os, Ana Rodr\'{\i}guez, Diana Stan, F\'elix del Teso, and Bruno Volzone.

\

\medskip

{\small
\bibliographystyle{amsplain}

}

\

\smallskip
{\footnotesize
\noindent {\sc Address}:  {Departamento de Matem\'aticas}\\[2pt]
   {Universidad Aut\'onoma de Madrid}\\[2pt]
   {28049 Madrid,  Spain}\\[2pt]
   {\sl e-mail:} {juanluis.vazquez@uam.es}
}

\vskip 1cm

2010 \textit{Mathematics Subject Classification.}
26A33, 
35K55, 
35K65, 
35S10. 

\medskip

\textit{Keywords and phrases.} Nonlinear diffusion, fractional Laplacian operator.

\end{document}